\documentclass[11pt]{amsart}
\usepackage{amssymb}
\usepackage{amscd}

\usepackage[all]{xy}
\newtheorem{theorem}{Theorem}[section]
\newtheorem{corollary}{Corollary}[section]
\newtheorem{lemma}{Lemma}[section]

\newtheorem{conjecture}{Conjecture}[section]

\newtheorem{remark}{Remark}[section]

\newtheorem{proposition}{Proposition}[section]
\newtheorem{example}{Example}[section]

\def \<{\langle}
\def \>{\rangle}

\def \a{\alpha }

\def \l{\lambda }

\newcommand{\bea}{\begin{eqnarray}}
\newcommand{\eea}{\end{eqnarray}}
\newcommand{\be}{\begin {equation}}
\newcommand{\ee}{\end{equation}}

\newcommand{\Z}{\Bbb Z}

\newcommand{\X}{\mathfrak{X}}
\newcommand{\W}{\mathcal W}
\newcommand{\mip}{\overline{M(1)}}

\newcommand{\Zp}{{\Bbb Z}_{>0} }

\newcommand{\N}{{\Bbb Z}_{\ge 0} }
\newcommand{\C}{\Bbb C}

\newcommand{\WW}{\boldsymbol{ \mathcal{W}}}

\newcommand{\la}{\langle}
\newcommand{\ra}{\rangle}

\newcommand{\tripletwp}{{\bf \mathcal{W}}_{2,p}}

\newcommand{\ds}{\displaystyle}

\begin{document}
 \title[On $\W$ algebra extensions of $(2,p)$ minimal models: $p > 3$ ]  {On $\W$-algebra extensions of $(2,p)$ minimal models: $p > 3$}
\author{Dra\v{z}en Adamovi\'c and Antun Milas}
\address{Department of Mathematics, University of Zagreb, Croatia}
\email{adamovic@math.hr}

\address{Department of Mathematics and Statistics,
University at Albany (SUNY), Albany, NY 12222}
\email{amilas@math.albany.edu}

\begin{abstract}
 This is a continuation of \cite{AdM-IMRN}, where, among other things, we classified irreducible representations
of the triplet vertex algebra $\mathcal{W}_{2,3}$. In this part we extend the classification to $\tripletwp$, for all odd $p>3$. We also determine the structure of the center of the Zhu algebra $A(\tripletwp)$ which implies the existence of a family of logarithmic modules having $L(0)$--nilpotent ranks $2$ and $3$.
A logarithmic version of Macdonald-Morris constant term identity plays a key role in the paper.
\end{abstract}

\maketitle

\section{Introduction}
Recently, there has been a stream of research on $\mathcal{W}$-algebras extensions of Virasoro minimal models from several different points of view \cite{AdM-lattice}, \cite{AdM-IMRN}, \cite{FGST-log}, \cite{FGST-log-kl}, \cite{GRW-1}. \cite{GRW-2}, \cite{PRZ}, \cite{R}, \cite{W}, etc.
As shown in \cite{FGST-log}, (see also \cite{AdM-IMRN}),  there is a remarkable vertex algebra $\WW_{q,p}$, $(p,q)=1$, an extension of the Virasoro vertex algebra
$L(c_{p,q},0)$ with central charge $c_{p,q}=1-\frac{6(p-q)^2}{pq}$. Studying the category of $\WW_{q,p}$-modules is interesting for several reasons.
On one hand, we expect to get new examples of $C_2$-cofinite, non self-dual conformal vertex algebras, which give rise to  finite tensor categories  \cite{Hu} \cite{HLZ} (although perhaps 
not necessarily rigid \cite{Miy}).  More interestingly,  the algebra $\WW_{q,p}$ is expected to be in Kazhdan-Lusztig duality with a certain quantum group $\mathcal{U}_{q,p}$ explicitly described in  \cite{FGST-log-kl}. In addition, representations of $\WW_{q,p}$ are rich in combinatorics and their
considerations should most likely lead to new parafermionic bases and $q$-series identities.

The main object of study in this paper are $W$-algebras $\WW_{2,p}$, $p \geq 2$.
In \cite{AdM-IMRN}, we start to investigate these $\mathcal{W}$--algebras in the framework of vertex algebras, where we proved their $C_2$--cofiniteness and irrationality. We also classified irreducible $\WW_{2,3}$-modules (this exactly correspond to the  $c=0$ triplet model from \cite{GRW-1} \cite{GRW-2}).  In the present paper we extend most of the results from \cite{AdM-IMRN} to every odd $p$.
An important role in our proofs is played by the {\em doublet} vertex superalgebra $\overline{V}_L$, with parity
decomposition
$$\overline{V}_L=\WW_{2,p} \oplus M,$$
where $M$ is a certain simple $\WW_{2,p}$-module.
Another key ingredient is the proof of Conjecture 10.1 from \cite{AdM-IMRN}.  
We obtained this via a logarithmic deformed version of Macdonald-Morris-Dyson constant term identities (see below).

Let us state the main results first.
\begin{theorem} \label{thm-nr}
The vertex algebra $\WW_{2,p}$ has precisely $4p+\frac{p-1}{2}$ (inequivalent) irreducible modules, explicitly constructed  in Section \ref{triplet-modules-sec}.
\end{theorem}

This result, together with explicit formulas of irreducible characters 
gives the following useful fact

\begin{theorem} \label{thm-mod} $SL(2,\mathbb{Z})$-closure of the space of irreducible $\WW_{2,p}$ characters is $\frac{15p-5}{2}$-dimensional
\end{theorem}

Having Theorem \ref{thm-nr} handy it is natural to ask for a complete structure of the Zhu algebra of $\WW_{2,p}$.
In our very recent work \cite{AdM-2010}, we introduced a new method for the description of Zhu's algebra for certain vertex algebras. This method was used to completely describe the structure of the Zhu algebra $A(\WW_{2,3})$. As an important consequence, we proved in \cite{AdM-2010} that $\WW_{2,3}$ admits a logarithmic module of $L(0)$--nilpotent rank $3$, conjectured previously in the physics literature.
This module is then used in the detailed analysis of the  $c=0$ triplet model in \cite{GRW-2}. In the present paper we apply the results from \cite{AdM-2010} for the Zhu algebra $A(\WW_{2,p})$. Although we cannot precisely describe $A(\WW_{2,p})$, we can still show

 \begin{theorem}  \label{thm-center} The center $Z(A( \WW_{2,p}))$ is isomorphic to
$$ {\C}[x] / \la f_{2,p} (x) \ra ,$$
where $f_{2,p}(x)$ is a certain polynomial of degree $\frac{15p-5}{2}$. 
\end{theorem}
Equality of dimensions in Theorem \ref{thm-mod} and \ref{thm-center} is far from accidental. As shown in \cite{FGST-log-kl} the center  of $\mathcal{U}_{2,p}$ is also $\frac{15p-5}{2}$-dimensional.
The previous result implies that $\WW_{2,p}$ admits $\frac{p-1}{2}$ (non-isomorphic) logarithmic modules of $L(0)$--nilpotent rank three.
Our forthcoming work \cite{AdM-2011} will provide explicit constructions of some logarithmic modules of $L(0)$-nilpotent rank three, including several
new modules.

Finally, we finish with a constant term identity, expressed as a residue identity. In a special case this identity is needed for purposes of proving Theorem \ref{thm-nr}.
\begin{theorem} \label{theorem-ct}
Let $k \geq 0$, and let also $p \geq 1$ be odd. Then
$${\rm Res}_{x_1,...,x_{2k+1}} \frac{\Delta(x_1,...,x_{2k+1})^p}{ (x_1 \cdots
x_{2k+1})^{(2k+1)p} } \prod_{i=1}^k {\rm ln} \left(1-\frac{x_{2i}}{x_{2i-1}} \right)
\prod_{i=1}^{2k+1} (1+x_i)^t$$
$$=\lambda_{k,p} \prod_{i=0}^{2k} {t+ \frac{pi}{2} \choose
(k+1)p-1},$$ where 
$$\Delta(x_1,...,x_{2k+1})=\ds{\prod_{1 \leq i < j \leq 2k+1}(x_i-x_j)},$$
and $\lambda_{k,p} \neq 0$ is a constant not
depending on $t$.
\end{theorem}
If we assume Conjecture \ref{log-dyson}, the constant $\lambda_{k,p}$ can be computed exactly.

\section{The vertex algebra $\tripletwp$ }

We assume the reader is familiar with vertex algebra theory as in say \cite{LL} of \cite{DL}.
In this section we shall consider the triplet vertex algebra $\WW_{2,p}$ introduced in  \cite{FGST-log} and \cite{AdM-IMRN} , a certain
subalgebra of the rank one lattice vertex algebra.

Assume that $p$ is an odd natural number, $p \ge 3$, and  let $$L = {\Z}
\alpha, \ \ \ \la \alpha , \alpha \ra = p.$$  As usual, we extend the scalars of $L$ and define an abelian algebra $\frak{h}=L \otimes \mathbb{C}$. We also write $\mathbb{C}[L]$
for the group algebra of $L$. Consider the usual affinization $\hat{\frak{h}}$ of $\frak{h}$, with brackets defined in the standard way by using the bilinear form on $L$. The Fock space of $\hat{\frak{h}}$ is denoted by 
$M(1)$. Then $$V_L \cong M(1) \otimes \mathbb{C}[L],$$ has a natural vertex superalgebra structure \cite{DL}, \cite{LL},
with vertex operator map
$$Y(u,x)=\sum_{n \in \mathbb{Z}} u_n x^{-n-1}.$$
Fix the Virasoro vector
$$\omega = \frac{1}{2 p} (\a (-1) ^2 + (p-2) \a (-2)) \in M(1) \subset V_L$$
and (screening) operators
$$Q = e ^{\a}_0, \qquad \widetilde{Q} = e ^{\frac{- 2 \a}{p} }_0, $$
$$G =\sum_{i =1 } ^{\infty} \frac{1}{i} e ^{\a}_{-i} e ^{\a}_i , \qquad G ^{tw}=\sum_{i =0 } ^{\infty} \tfrac{1}{i+ 1/2} e ^{\a}_{-i-1/2} e ^{\a}_{i+ 1/2} .$$
The action of $G$ (resp. $G ^{tw}$)
is well-defined on any $V_L$-module (resp. $\Z_2$--twisted $V_L$--module). For details see \cite{AdM-IMRN}.

Although $\WW_{2,p}$ was originally defined as
the intersection of two screening  operators \cite{FGST-log}, we showed in \cite{AdM-IMRN} that it can be realized as a subalgebra of $V_L$ generated by $\omega$ and three primary vectors
$$ F = Q e ^{-3 \a}, \quad H = G F, \quad E = G ^2 F.$$
The algebra $\WW_{2,p}$ is of course ${\Z}$--graded whose charge zero component is the singlet vertex algebra $\mip$ generated by $\omega$ and $H$.
Let
$ A(\mip) = \mip / O(\mip)$ and $A(\WW_{2,p}) = \WW_{2,p} / O(\WW_{2,p})$ be the associated Zhu algebras of $\mip$ and $\WW_{2,p}$, respectively  (see \cite{AdM-IMRN} for details).  

We recall further results from \cite{AdM-IMRN} on the structure of the Zhu algebra $A(\mip)$.
 Let
$$h_{r,s}=\frac{(pr-2s)^2-(p-2)^2}{8p}.$$

  \begin{theorem} \label{zhu-singlet-2-p}
The Zhu algebra $A(\overline{M(1)})$
is isomorphic to the commutative associative algebra
$$A(\overline{M(1)}) \cong \frac{{\C}[x,y]}{ \langle P(x,y) \rangle }$$
where
$$ P(x,y) = y ^2 - C_p  \prod_{i=1} ^{{2p-1}} \left( x- h_{1,i}\right) ^2 \prod_{i = 1} ^{2p-1} \left( x- h_{2,i} \right)\qquad
(C_p \ne 0).$$
(Here $x$ corresponds to $[\omega]$ and $y$ to $[H]$.)
\end{theorem}

\section{Some modules for the doublet vertex superalgebra $\overline{V_L}$}

As in \cite{AdM-IMRN}, we consider the vertex superalgebra $\overline{V_L}$  generated by $\omega$, 
$$ a ^{-} = Q e ^{-2\a} \quad \mbox{and} \quad a ^+ = G a ^- . $$
In this section we shall describe some twisted and untwisted modules for $\overline{V_L}$. Some proofs in this section are analogous to the proofs of similar results in Section 5 of  \cite{AdM-IMRN}, so we omit some details for brevity. We should also say that these results are in agreement with the structural results obtained in \cite{FGST-log}. Main difference in our approach is in the fact that some  complicated screening operators constructed in \cite{TK} are now replaced by exponents of screening operators $G$ and $G ^{tw}$.

The proof of the following theorem is completely analogous to that of Theorem 4.1 of \cite{AdM-IMRN}, so we only indicate the main steps of the proof.
\begin{theorem} \label{fr-p-gen} Let $1 \le k \le p$.
The space of intertwining operators
\bea \label{03p2} && I { L(c_{2,p},h) \choose
L(c_{2,p},h_{5,1}) \ \ L(c_{2,p},h_{n,k} ) }\eea
is nontrivial only if
\be \label{03p2-fr} h \in \left\{ h_{n-4,k}, h_{n-2,k},
h_{n,k}, h_{n+2,k}, h_{n+4,k}
%
%
 \right\}. \ee

\end{theorem}
\begin{proof} (Sketch) As in \cite{AdM-IMRN}, we have a singular vector $v_{sing} \in M(c_{2,p},h_{5,1})$ of (relative) degree five, which generates a submodule 
$M_1 \subset M(c_{2,p},h_{5,1})$. Again, as in \cite{AdM-IMRN}, we use the Frenkel-Zhu's formula to show that if space
$I { L(c_{2,p},h) \choose M(c_{2,p},h_{5,1})/M_1 \ \ L(c_{2,p},h_{n,k} ) }$ is nontrivial then
$h$ must be in the given range.
But then the same holds for any quotient of $M(c_{2,p},h_{5,1})/M_1$ and in particular for $L(c_{2,p},h_{5,1})$ . The proof follows.
\end{proof}

We will need another result from \cite{AdM-IMRN}
\begin{proposition} \label{CT} Assume that $v_{\l}$ is a lowest weight vector in $M(1)$--module $M(1, \l)$ such that $\alpha(0)v_{\l}=\lambda(\alpha)v_{\l}$. Let $ t = \l (\a)$. Then we have:
$$ H(0) v_{\l} = D_p { t + p \choose 2p-1} {t \choose 2p -1} { t + p/2
\choose 2p -1} v_{\l} \quad (D_p \ne 0).$$
\end{proposition}

Define now the following cyclic $\overline{V_L}$--modules:

$$ A_p( h_{1,2p-k} ) := \overline{V_L} . Q e ^{ \tfrac{k-1}{p} \a -\a}, \quad A_p (h_{2,k} ) := \overline{V_L} . e ^{-\tfrac{\a}{2} + \tfrac{k-1}{p} \a} \quad (k=1, \dots,p). $$

Then $A_p( h_{1,2p-k} )$ is an  untwisted and $A_p (h_{2,k} )$ is  a $\Z_2$--twisted $\overline{V_L}$--module.
Let us denote by $Y(\cdot,z)$ and $Y^{tw}(\cdot,z)$ the associated vertex operators.

The key point for the description of these modules are in the  following lemmas which gives non-triviality of certain singular vectors.
The proof of the first lemma is  identical to that of Lemma 5.1 of \cite{AdM-IMRN}.
\begin{lemma} \label{non-triv-gen1} Assume that $1 \le k \le p$. We have:
\item[(i)] $Y(a ^-, x) Q e ^{ \tfrac{k-1}{p} - (n+1) \a } \in W ((x)), $
where $$W = U  (Vir). Q e^{ \tfrac{k-1}{p} \a  -(n+2) \a} \cong L ^{Vir}(c_{2,p},
h_{2n+5,k}).$$

\item[(ii)] $ G^{n} Q e^{\tfrac{k-1}{p}\a -(n+1) \a} \ne 0 \quad  \mbox{for} \ n
\in {\N}$,  $ G ^j  Q e^{\tfrac{k-1}{p}\a -(n+1) \a} = 0$  for $ j > n$.
\item[(iii)] $ G^{n} Q e^{\tfrac{k-1}{p}\a-(n+1) \a} \in A_p( h_{1,2p-k} )$.
\end{lemma}

\begin{lemma} \label{non-triv-gen2} Assume that $1 \le k \le p$. We have:
\item[(i)] $Y^{tw} (a ^-, x) e ^{ -\tfrac{\a}{2} + \tfrac{k-1}{p} \a - n \a } \in W ((x)), $
where $$W = U  (Vir). e ^{ -\tfrac{\a}{2} + \tfrac{k-1}{p}\a - (n+1) \a } \cong L ^{Vir}(c_{2,p},
h_{2n+4,k}).$$

\item[(ii)] $ (G^{tw} ) ^n   e^{-\tfrac{\a}{2} + \tfrac{k-1}{p} \a -n  \a} \ne 0 \quad  \mbox{for} \ n
\in {\N}$,   $ (G^{tw} ) ^j   e^{-\tfrac{\a}{2} + \tfrac{k-1}{p} \a -n  \a} =  0 \quad  \mbox{for}  \ j > n$.
\item[(iii)] $ (G^{tw} ) ^n   e^{-\tfrac{\a}{2} + \tfrac{k-1}{p}\a -n  \a} \in A_p (h_{2,k} )$.
\end{lemma}
\begin{proof}
First we shall prove assertion (i).
Since $V_L$ is a simple vertex superalgebra and $V_{L - \a /p + \tfrac{k-1}{p} \a}$ is its simple twisted module we conclude that
$$Y^{tw} (a ^-, x) e ^{ -\tfrac{\a}{2} + \tfrac{k-1}{p} \a - n \a } \ne 0.$$ Therefore, we can  find $j_0 \in \tfrac{1}{2} + \Z$, such that
$$a ^- _{j_0}  e ^{-\tfrac{\a}{2} + \tfrac{k-1}{p}\a - n \a} \ne 0,  \quad a ^- _{j}  e ^{-\tfrac{\a}{2} + \tfrac{k-1}{p}\a - n \a} = 0 \quad \mbox{for} \ j >j_0.$$

By using fusion rules from Theorem \ref{fr-p-gen} we conclude that $a ^- _{j_0}  e ^{-\tfrac{\a}{2} + \tfrac{k-1}{p}\a - n \a}$ must be a singular vector in $
M(1) \otimes e ^{ -\tfrac{\a}{2} + \tfrac{k-1}{p}\a - (n +1)\a}$ of conformal weight $h_{2n+4,k}$ (there are no singular vectors of weight $h_{2n+6,k}$ in this Fock space). Therefore
$j_0 = - (n- 3/2)
p + k -3$ and
$$ a ^- _{j_0}  e ^{-\tfrac{\a}{2} + \tfrac{k-1}{p}\a - n \a} = \mu_n   e ^{ -\tfrac{\a}{2} + \tfrac{k-1}{p}\a - (n +1)\a} \quad
( \mu_n \ne 0), $$
$$ a^- _ j e ^{-\tfrac{\a}{2} + \tfrac{k-1}{p}\a - n \a}  \in W \quad \mbox{for} \ j \le j_0.  $$
 In this way we have proved assertion (i).

 We shall prove the assertions (ii) and (iii) by induction on $n \in
{\Zp}$.

Since $H(0) e ^{-\tfrac{a}{2} + \tfrac{k-1}{p} \a - \a} \ne 0$ (see Proposition \ref{CT}), we conclude that
$$G ^{tw} e ^{-\tfrac{a}{2} + \tfrac{k-1}{p} \a - \a} = \mu_0   a^+ _ { 3 p /2 +k-3} e ^{-\tfrac{\a}{2} + \tfrac{k-1}{p}\a } \ne 0. $$
So the assertion holds for $n =1$.

Assume now that assertions (ii)-(iii)  hold for certain $n \in
{\Zp}$. Since $V_{L} $ is a simple vertex operator algebra we have
that
$$ Y(a ^{+},z) (G^{tw} )^{n}  e ^{-\tfrac{\a}{2} + \tfrac{k-1}{p}\a - n \a} \ne 0,$$
(for the proof see \cite{LL}).
So there is $k_0 \in \tfrac{1}{2} + {\Z}$ such that
$$a ^+ _{k_0} (G^{tw} )^{n}  e ^{-\tfrac{\a}{2} + \tfrac{k-1}{p}\a - n \a}   \ne 0 \quad \mbox{and} \quad
 a ^+ _{j} (G^{tw} )^{n}  e ^{-\tfrac{\a}{2} + \tfrac{k-1}{p}\a - n \a}  = 0 \ \mbox{for} \ j > k_0.$$
Since
$$ a ^+ _{k_0} (G^{tw} )^{n}  e ^{-\tfrac{\a}{2} + \tfrac{k-1}{p}\a - n \a} = \nu_2 (G^{tw} )^{n+1} (a ^{-} _{k_0}
 e ^{-\tfrac{\a}{2} + \tfrac{k-1}{p}\a - (n+1) \a} ), $$
for certain non-zero constant ${\nu}_2$, then by using assertion
(i) and the fact that $ G ^{tw} $ is a screening operator we conclude that
$$ a ^+ _{k_0} (G^{tw} )^{n}  e ^{-\tfrac{\a}{2} + \tfrac{k-1}{p}\a - n \a}   \in U(Vir) (G^{tw} )^{n+1}  e ^{-\tfrac{\a}{2} + \tfrac{k-1}{p}\a - (n+1) \a}. $$
Therefore
$(G^{tw} )^{n+1}  e ^{-\tfrac{\a}{2} + \tfrac{k-1}{p}\a - (n+1) \a} \ne 0$ and
$$ (G^{tw} )^{n+1}  e ^{-\tfrac{\a}{2} + \tfrac{k-1}{p}\a - (n+1) \a}  = \frac{1}{{\nu}_2 {\mu_n}} a ^+ _{j_0} (G^{tw} )^{n}  e ^{-\tfrac{\a}{2} + \tfrac{k-1}{p}\a - n \a} \in  A_p (h_{2,k} ).$$
The proof follows.

\end{proof}

By using fusion rules from Theorem \ref{fr-p-gen}, Lemmas \ref{non-triv-gen1} and \ref{non-triv-gen2}  and similar methods as in \cite{AdM-IMRN} we can describe these modules:

\begin{theorem} \label{decompos-a}
\item[(1)] $A_p( h_{1,2p-k} )$ is a completely reducible module for the Virasoro algebra generated by singular vectors:
$$ v_{k,n}  ^{(1)} =G ^j Q  e ^{ \tfrac{k-1}{p} \a - (n+1) \a}, \ n \ge 0, \ 0 \le j \le n. $$
We have the following decomposition:
$$ A_p( h_{1,2p-k} )= \bigoplus _{n=0} ^{\infty} (n+1) L (c_{2,p}, h_{2n+3,k} ). $$
\item[(2)] $A_p( h_{2,k} )$ is a completely reducible module for the Virasoro algebra generated by singular vectors:
$$ v_{k,n}  ^{(2)} =(G^{tw}) ^j   e ^{ -\tfrac{\a}{2} + \tfrac{k-1}{p} \a - n \a}, \ n \ge 0, \ 0 \le j \le n. $$
We have the following decomposition:
$$ A_p( h_{2,k} )= \bigoplus _{n=0} ^{\infty} (n+1) L (c_{2,p}, h_{2n+2,k} ). $$
\end{theorem}

We have the following useful description of  the structure of some $\overline{V_L}$--modules constructed above.

\begin{proposition} We have:
\bea  && A_p (h_{1,p}) = \mbox{Ker}_{ V_{L + \tfrac{p-1}{p} \a} }   \ Q , \\
&& A_p (h_{2,1}) = \mbox{Ker}_{V_{L -  \tfrac{\a}{2}  } } \widetilde{Q}, \\
&& A_p (h_{2,p}) =  {V_{L -\tfrac{\a} {2} +    \tfrac{p-1} {p} \a } } . \eea
\end{proposition}

\section{Some $\tripletwp$-modules}
\label{triplet-modules-sec}

In this part we introduce some $\WW_{2,p}$-modules and describe their structure as $Vir$-module.
In a special case we obtain Proposition 5.4 in \cite{AdM-IMRN}, a decomposition of $\WW_{2,p}$.

Let $V_L$ be as before. We also let $D=2 \mathbb{Z} \alpha$, so that $V_D$ is a vertex subalgebra of $V_L$.
There are $4p$ irreducible $V_L$-modules, corresponding to $4p$ cosets in the quotient $L^\circ/L$, which we
conveniently describe here.

For $k=1, \dots ,p$ consider:

$$V_{D+(k-1)\alpha/p}, \ \ V_{D+(k-1)\alpha/p-\alpha}, \ \ V_{D - \a /2 + (k-1)\alpha/p}, \ \ V_{D - \a /2 + (k-1) \a /p - \a  }.$$
By restriction, these are also modules for the triplet $\WW_{2,p}$, usually called "Verma modules".

The singular vectors in $V_{D+\gamma}$ generate ${\rm Soc}_{Vir}(V_{D+\gamma})$. We let
\bea
&& \X^+_{1,k}:={\rm Soc}_{Vir} (V_{D+(k-1)\alpha/p}), \nonumber \\
&& \X^+_{2,k}:={\rm Soc}_{Vir} (V_{D- \alpha/2 + (k-1)\alpha/p}), \nonumber \\
&& \X^-_{1,k}:={\rm Soc}_{Vir} (V_{D+(k-1)\alpha/p-\alpha}),\nonumber \\
&& \X^-_{2,k}:={\rm Soc}_{Vir} (V_{D- \a /2 +(k-1)\alpha/p - \a }), \nonumber
\eea

Our notation clearly follows parametrization from \cite{FGST-log}.

We also let $M(c,h)$ denote the Virasoro Verma module of lowest weight $h$ and central charge $c$, and denote
by $L(c,h)$ the corresponding irreducible quotient.

By combining results from \cite{FGST-log} and action of screening operators $G$ and $G^{tw}$ from  Lemmas \ref{non-triv-gen1} and \ref{non-triv-gen2} (see also \cite{AdM-IMRN})   we get the following result:

\begin{theorem} \label{decompos}   As a $Vir$-modules  $\X^{\pm} _{s,r}$ are generated by the family of singular vectors $\mbox{Sing} ^{\pm} _{s,k}$, where
 $s\in \{1,2\}$, $k \in \{1, \dots, p \}$ and
 we have:
\bea
&&
\mbox{Sing} ^{+} _{1,k} =\{ G ^ j Q e ^{ \tfrac{k-1}{p} \a - (2n+1) \a} \ \vert \ n \in {\N}, \ 0 \le j \le 2n \},\nonumber \\
&&
\mbox{Sing} ^{-} _{1,k} =\{ G ^ j Q e ^{ \tfrac{k-1}{p} \a - 2n \a} \ \vert \ n \in {\Zp}, \ 0 \le j \le 2n-1 \}, \nonumber \\
&&
\mbox{Sing} ^{+} _{2,k} =\{ (G^{tw}) ^ j  e ^{ -\tfrac{1}{2} \a + \tfrac{k-1}{p} \a - 2n \a} \ \vert \ n \in {\N}, \ 0 \le j \le 2n \},\nonumber \\
&&\mbox{Sing} ^{-} _{2,k} =\{ (G^{tw}) ^ j  e ^{ -\tfrac{1}{2} \a +  \tfrac{k-1}{p} \a - (2 n - 1) \a} \ \vert \ n \in {\Zp}, \ 0 \le j \le 2n-1 \}.\nonumber
\eea
We have the following decompositions:
\bea
&& \X^+_{1,k}=\oplus_{n=0}^\infty (2n+1)L(c_{2,p}, h_{4n +3, k}) \nonumber \\
&& \X^+_{2,k}=\oplus_{n=0}^\infty (2n+1)L(c_{2,p}, h_{4n +2,k})\nonumber \\
&& \X^-_{1,k}=\oplus_{n=1}^\infty (2n)L(c_{2,p}, h_{4n +1, k}) \nonumber \\
&& \X^-_{2,k}=\oplus_{n=1}^\infty (2n)L(c_{2,p}, h_{4n,k}). \nonumber
\eea
\end{theorem}

In order to see that $ \X ^{\pm} _{s,r}$ has the structure of $\WW_{2,p}$--modules we shall use results from previous section and construction of $\overline{V_L}$--modules.
The vertex superalgebra $\overline{V_L}$ has a canonical $\Z_2$--automorphism $\sigma$ and the fixed point subalgebra is exactly the triplet vertex algebra $\WW_{2,p}$. Moreover every (twisted) $\overline{V_L}$--module from previous section is a $\Z_2$--graded and decomposes into direct sum of two ordinary  $\WW_{2,p}$--modules.
We have:
\bea
&& A_p(h_{1,2p-k}) = W_{p} (h_{1,2p-k}) \oplus W_{p} ( h_{1,3p-k}), \nonumber \\
&& A_p(h_{2,k} )= W_{p} (h_{2,k}) \oplus W_{p} ( h_{2,3p-k}). \nonumber
\eea
In this way we have obtained four families of ordinary $\WW_{2,p}$--modules.
By using previous results on the structure of $\overline{V_L}$--modules one can easily obtains that the constructed $\WW_{2,p}$--modules are cyclic and
\bea
&& W_{p} (h_{1,2p-k}) = \WW_{2,p} . Q e ^{ \tfrac{k-1}{p} \a - \a}, \label{mod-1} \\
&& W_{p}(h_{1, 3p-k} ) = \WW_{2,p} .  Q e ^{ \tfrac{k-1}{p} \a - 2 \a}, \label{mod-2}\\
&& W_{p} (h_{2,k}) = \WW_{2,p} . e ^{ -\tfrac{ \a}{2} + \tfrac{k-1}{p} \a }, \label{mod-3} \\
&& W_{p}(h_{2, 3p-k} ) = \WW_{2,p} .   e ^{ - \tfrac{\a}{2} + \tfrac{k-1}{p} \a -  \a}. \label{mod-4}
\eea
Now by using Theorem \ref{decompos-a} and Theorem \ref{decompos}  we get the following result.
\begin{corollary}
$\X ^{\pm}_{1,k}$, $\X ^{\pm} _{2,k}$  are $\WW_{2,p}$--modules and
$$ \X ^{+}_{1,k} = W_{p} (h_{1,2p-k}), \ \X^{-} _{1,k} =  W_{p}(h_{1, 3p-k} ), \ \X ^{+}_{2,k} = W_{p} (h_{2,k}), \ \X ^{-} _{2,k} = W_{p}(h_{2, 3p-k} ). $$
\end{corollary}

\vskip 5mm

We also let

\begin{equation} \label{k-def}
\mathcal{K}^+_{1,r}={V}(c_{2,p},\frac{(p-2r)^2-(p-2)^2}{8p})+\X^+_{1,r},
\end{equation}
where ${V}(c_{2,p},\frac{(p-2r)^2-(p-2)^2}{8p})$ is the $Vir$-module generated by $e^{(r-1)\alpha/p}$.

More precisely,  ${V}(c_{2,p},\frac{(p-2r)^2-(p-2)^2}{8p})$ is a quotient of the Verma module
${M}(c_{2,p},\frac{(p-2r)^2-(p-2)^2}{8p})$, and it combines in an exact sequence
$$0 \rightarrow L(c_{2,p}, \frac{(3p-2r)^2-(p-2)^2}{8p})  \rightarrow  {V}(c_{2,p},\frac{(p-2r)^2-(p-2)^2}{8p})$$
$$ \rightarrow {L}(c_{2,p},\frac{(p-2r)^2-(p-2)^2}{8p}) \rightarrow 0.$$
If we restrict ourself in the range $1 \leq r \leq \frac{p-1}{2}$, this way we obtain $(2,p)$-minimal models.

\begin{proposition} \label{podmodul}
\item[(1)] We have that
$\mathcal{K}^+_{1,k} = \WW_{2,p}.  e ^{ \tfrac{k-1}{p} \a}. $

 \item[(2)] $ W_{p} (h_{1,2p - k})$ is an maximal submodule of  $\mathcal{K}^+_{1,k}$ and
\bea && W_p (h_{1,k}) := L(c_{2,p}, h_{1,k} ) = \frac{ \mathcal{K}^+_{1,k}     }{ W_{p} (h_{1,2p-k})   } \quad (k=1, \cdots, \frac{p-1}{2}). \label{min-mod} \eea
is isomorphic to the Virasoro minimal model with lowest  conformal weight  $h_{1,k}$.
\end{proposition}
\begin{proof}
Assertion (1) follows from the fact that $ Q e ^{\tfrac{k-1}{p} \a - \a}$ is a singular vector in ${V}(c_{2,p},  h_{1,k})$ which generates its maximal submodule.

Let $X \in \{ E, F, H \}$. Set $X(n) = X_{ 6 p - 4 + n}$.
First we notice that
$$ L(n+1)  Q e^{\frac{k-1}{p} \a -\a } = X (n)  Q e^{\frac{k-1}{p} \a -\a } = 0 \qquad ( n \in \N). $$
The first relation is clear. To see that $F(n)Q e^{\frac{k-1}{p} \a -\a }=0$, for $n \in \N$, it is sufficient
to observe $Q^2 e^{-3 \alpha}=0$  and $Q^2 e^{(k-1)/p\alpha-\alpha}=0$. For $X=H$ and $X=F$,  apply
Lemma \ref{non-triv-gen1} (ii).

Therefore $e^{\frac{k-1}{p} \a  } \notin W_{p} (h_{1,2p - k})$ and $W_{p} (h_{1,2p - k})$ is a maximal  submodule of $\mathcal{K}^+_{1,k}$.
The proof follows.
\end{proof}

\section{Classification of irreducible $\tripletwp$-modules}

In   \cite{AdM-IMRN} we  showed that classification of irreducible $\WW_{2,p}$--modules is related to a
constant term identities. In fact Conjecture 10.1. from \cite{AdM-IMRN} essentially gives classification of irreducible modules. In Section \ref{identitet} below we shall prove that this conjecture holds. Now we shall briefly discuss the classification of irreducible modules.

 The results from Section 10 of \cite{AdM-IMRN} give the following relation in
$A( \WW_{2,p})$:
\begin{theorem} \label{klasif-general}In the Zhu algebra $A( \WW_{2,p})$ we have:

$$f_{2,p}([\omega] ) = 0,$$
where
\bea f_{2,p}(x) &=&\left( \prod _{i =1} ^{3p-1} (x - h_{1,i} )  \right) \left(  \prod_{i =1} ^{
 \tfrac{3p-1}{2}
} (x - h_{1,2p-i} )\right) \left(\prod_{i =1} ^{3p-1} (x - h_{2,i} )\right )  \nonumber \\
&=& \left(\prod_{i =1} ^{ \tfrac{p-1}{2}} (x- h_{1,i}) \right)^3  \left(\prod_{i = p } ^{ 2p -1} (x- h_{1,i}) \right)^2 \left(\prod_{i =1} ^{p-1} (x - h_{2,i} )\right ) ^2   \nonumber \\
&& \cdot (x-h_{2,p}) \left(\prod_{i=2p} ^{ 3 p -1} (x- h_{1,i}) \right) \left(\prod_{i =2 p } ^{3 p-1} (x - h_{2,i} )\right ).  \eea
\end{theorem}

Following approaches in \cite{AdM-triplet} and \cite{AdM-IMRN} we can now classify all irreducible $\WW_{2,p}$. We shall omit some details. Let

 $$S_{2,p} =   \{ h_{1,i}, h_{1,j}, h_{2,k}, h_{2,l}  \vert \ 1 \le i \le \tfrac{p-1}{2},  \ p \le j \le 3p-1, \ 1 \le k \le p, \ 2p \le l \le 3p-1  \}. $$

\begin{theorem} \label{constr-irr}  We have:
\item[(1)]
$ W_{p} (h_{1,2p - k})$  is an irreducible   $\N$--graded $\WW_{2,p}$--module with lowest weight $h_{1,2p -k}$. Its top component is $1$-dimensional irreducible module for the Zhu algebra $A(\WW_{2,p})$ spanned by $Q e ^{ \tfrac{k-1}{p} \a-\a}$.

\item[(2)]
$ W_{p} (h_{1,3p - k})$  is an irreducible   $\N$--graded $\WW_{2,p}$--module with lowest weight $h_{1,3p -k}$. Its top component is $2$-dimensional irreducible module for the Zhu algebra $A(\WW_{2,p})$ spanned by $Q e ^{ \tfrac{k-1}{p} \a - 2 \a}$ and $ G Q e ^{ \tfrac{k-1}{p} \a - 2 \a}$.

 \item[(3)]
$ W_{p} (h_{2, k})$  is an irreducible   $\N$--graded $\WW_{2,p}$--module with lowest weight $h_{2,k}$. Its top component is $1$-dimensional irreducible module for the Zhu algebra $A(\WW_{2,p})$ spanned by $e ^{ - \tfrac{\a}{2}+ \tfrac{k-1}{p} \a  }$

 \item[(4)]
$ W_{p} (h_{2, 3 p -k})$  is an irreducible   $\N$--graded $\WW_{2,p}$--module with lowest weight $h_{2,3p -k}$. Its top component is $2$-dimensional irreducible module for the Zhu algebra $A(\WW_{2,p})$ spanned by $ e ^{ - \tfrac{\a}{2}+ \tfrac{k-1}{p} \a -  \a}$ and $ G ^{tw}  e ^{ - \tfrac{\a}{2}+ \tfrac{k-1}{p} \a -  \a}$.
\end{theorem}
\begin{proof}
Let us prove assertion (1). In Proposition \ref{podmodul} we proved that $ W_{p} (h_{1,2p - k})$ is $\N$--graded and that its top component is spanned by $Q e ^{ \tfrac{k-1}{p} \a - \a}$.
Next we notice that the  set conformal weights of  singular vectors    appearing in the decomposition of $ W_{p} (h_{1,2p - k})$
is   $ \{  h_{ 4 n +3, k }  \ \vert    n \ge 0 \} $  and that
\bea && h_{4 n +3, k } \notin S_{2,p} \qquad \mbox{for} \ n \ge 1. \label{uv-4} \eea
Assume now that $N \subset
W_{p} (h_{1,2p - k})$ is a non-trivial submodule. Then $N$ is a
${\N}$--graded, and  Theorem \ref{klasif-general} gives that the top component $N(0)$ must have conformal weight
$h \in S_{2,p}$. On the other hand every non-trivial vector from the top component is a singular vector for the Virasoro algebra and therefore  $h = h_{4n+3,}$ for certain $n \ge 1$. This contradicts  (\ref{uv-4}). The proof of other assertions are similar.
\end{proof}

By using Theorems \ref{klasif-general} and  \ref{constr-irr}  and the same proof as in \cite{AdM-triplet} and \cite{AdM-IMRN} we get:

 \begin{theorem}
 The set
 $$ \{ W_p (h), \quad h \in S_{2,p} \}$$
 provides, up to isomorphism, all irreducible $\WW_{2,p}$--modules.
 \end{theorem}

As in \cite{AdM-2010} we have a natural homomorphism
\bea  \Phi : && A(\mip)  \rightarrow  A(\WW_{2,p}) \nonumber \\
&& a + O(\mip)    \mapsto  a + O(\WW_{2,p})  \qquad (a \in \mip) . \nonumber \eea

By using the description of the Zhu algebra $A(\mip)$ from \cite{AdM-IMRN} and the methods from \cite{AdM-2010} we get obtain the following result:

\begin{proposition} \label{ker-c2}
$\mbox{Ker}(\phi)$ is contained in the following ideal in $A(\mip)$
$$ A(\mip) . \{ p([\omega]) * [H], f_{2,p}([\omega]) \} = \{ A([\omega]) p([\omega]) [H] + B([\omega]) f_{2,p}([\omega]), \ \ A, B \in {\C}[x] \}, $$
where
 \bea  p(x) =  \prod_{i=2p} ^{ 3 p -1} (x- h_{1,i}) \prod_{i=2p} ^{ 3 p -1} (x- h_{2,i}). \nonumber \eea
\end{proposition}

\begin{theorem} \label{klasif-general-center}  The center of  the Zhu algebra $A( \WW_{2,p})$ is isomorphic to
$$ {\C}[x] / \la f_{2,p} (x) \ra $$
\end{theorem}
\begin{proof}
By using results from \cite{AdM-IMRN} one can easily see that the center of $A(\WW_{2,p})$ is isomorphic to the subalgebra generated by $[\omega]$.
By using Proposition \ref{ker-c2} we see that
$$ f \in {\C}[x], \ f([\omega]) \in \mbox{Ker} (\Phi) \quad \implies  f_{2,p} \vert f. $$
Since $G ^2  Q e ^{-5\a} \in O(\WW_{2,p})$ and
$$ [ G ^2  Q e ^{-5\a} ] = \nu f_{2,p}([\omega]) \quad \mbox{in} \ A(\mip) \qquad (\nu \ne 0)$$
we conclude that
 $f_{2,p} ([\omega]) \in \mbox{Ker}(\phi)$. The proof follows.
\end{proof}

The structure of the center easily implies that
 $A(\WW_{2,p})$ contains two and three dimensional indecomposable modules.

 \begin{proposition} The Zhu algebra $A(\WW_{2,p})$ has  $2$--dimensional idecomposable modules $U^{(2)} _{h_{1,i}}$ on which $[\omega]$ acts as
$$ \left(
   \begin{array}{cc}
     h_{1,i} & 1 \\
     0 & h_{1,i} \\
   \end{array}
 \right) $$
  where  $p \le i  \le 2p-1$  and $3$--dimensional indecomposable modules $U^{(3)} _{h_{1,i}}$ on which $[\omega]$ acts as
 $$
 \left(
   \begin{array}{ccc}
     h_{1,i} & 1 & 0 \\
     0 & h_{1,i} & 1 \\
     0 & 0 & h_{1,i} \\
   \end{array}
 \right) $$
 where  $1 \le i  \le (p-1) /2.$
\end{proposition}

 By applying the theory of Zhu's algebras, we get the existence of   logarithmic, $\N$--graded   $\WW_{2,3}$--modules $R_h ^{(2)}$ and $R_h ^{(3)}$ whose top components are isomorphic to $U_h ^{(2)}$ and $U_h ^{(3)}$ respectively.
 Therefore we have  proved:

\begin{corollary} The vertex algebra $\WW_{2,p}$ admits $(p-1)/2$ logarithmic modules
of $L(0)$ nilpotent rank $3$.
\end{corollary}

\begin{remark}
It is important problem to construct all logarithmic modules for triplet vertex algebras. In this section we proved the existence of all logarithmic modules which can be detected by using Zhu algebra theory. In \cite{AdM-2011} we  present an explicit construction of certain logarithmic modules of $L(0)$ nilpotent rank 3.
\end{remark}

\section{ The structure of some irreducible $\mip$--modules}

In this section we shall investigate some modules for the singlet vertex algebra $\mip$. Since $\mip$ is subalgebra of the free boson vertex algebra $M(1)$, it is interesting to find   $M(1)$--modules $M(1) \otimes e ^{\l}$ which  will be irreducible when we restrict them to $\mip$. The similar modules in the case of $(1,p)$--modules were investigated in  \cite{A-2003} and \cite{AdM-2007}.

\begin{proposition}
\item[(1)] For every $m \in {\Z}$ $M(1) \otimes e ^{ \tfrac{ p-2} {2p}  \a  + m \a}$ is an irreducible $\mip$--modules.

\item[(2)] $M(1) \otimes e ^{ \tfrac{ p-2} {2p}  \a }$ is an irreducible self-dual $\mip$--module. As a Virasoro module,

\bea   M(1) \otimes e ^{ \tfrac{ p-2} {2p}  \a } &\cong&
\bigoplus_{n=0} ^{\infty}  U(Vir). (G ^{tw}) ^n e ^{ \tfrac{p-2}{2p} \a - 2n \a} \nonumber \\
& \cong &\bigoplus_{n=0} ^{\infty} L(c_{2,p}, h_{ 4n +2,p} ). \nonumber \eea
\end{proposition}
\begin{proof}
First we notice that
$$\X _{2,p} ^{+} = V_{ D + \tfrac{ p-2} {2p}  \a} \quad \mbox{and} \quad  \X _{2,p} ^{-} = V_{ D + \tfrac{  3 p-2 }{2p} \a} $$
are irreducible $\WW_{2,p}$--modules. $\X_{2,p} ^{\pm}$ are naturally $\Z$--graded and it is easy to see that every graded component must be an irreducible module for $\mip$. The proof follows.
\end{proof}

\section{Logarithmic Macdonald-Morris-Dyson identities} \label{identitet}

In \cite{AdM-IMRN}, we conjectured certain constant term identities (or residue identity) involving 
\be \label{ctt}
{\rm Res}_{x_1,...,x_{2n+1}} \frac{1}{(x_1 \cdots x_{2n+1})^i} \prod_{i<j} (x_i-x_j)^p  \prod_{i=1}^{2n+1} (1+x_i)^t  \prod_{i=1}^n {\rm ln}\left(1-\frac{x_{2i}}{x_{2i-1}} \right),
\ee
where $p$ is positive and {\em odd}, $i \in \mathbb{Z}$, and ${\rm ln}\left(1-\frac{x}{y} \right)=-\sum_{m=1}^\infty \frac{x^m}{m y^m}$.
Logarithmic factors are really needed here (if absent, the residue  (\ref{ctt}) would be trivially zero!).
Related identities without logarithms factors have been studied in the literature for a while.
Let us first notice that the rational function part in (\ref{ctt}) resembles rational function appearing in celebrated  Dyson-Macdonald's constant term identities for root systems of  type $A$  (cf. \cite{FW})
$${\rm CT}_{x_1,...,x_n} \prod_{1 \leq i \neq j \leq n} \left(1-\frac{x_i}{x_j}\right)^m=\frac{(mn)!}{m!^n},$$
 which holds for every $m \in \mathbb{N}$. But (\ref{ctt}) involves an additional parameter 
 $t$,  so it should viewed as logarithmic deformation not of Dyson's  but of Morris constant term identity \cite{M} 
\bea \label{morris}
{\rm CT}_{x_1,...,x_n}  \prod_{i=1}^n (1+x_i)^a \left(1+\frac{1}{x_i} \right)^b \prod_{1 \leq i \neq j \leq n} \left(1-\frac{x_i}{x_j} \right)^m=\prod_{i=0}^n \frac{(a+b+mi)!(m(i+1))!}{(a+mi)! (b+mi)! m!}, \nonumber
\eea
which is valid for every $a,b,m \in \mathbb{N}$. 
Although there  seems to be no direct link between evaluation of (\ref{morris}) and (\ref{ctt}) for the above reasons we refer to evaluation of (\ref{ctt}) as {\em logarithmic Macdonald-Morris} identity.

In what follows, we shall evaluate (\ref{ctt}) up to a nonzero scalar (which depends on $p$ but not on $t$).

For start, we denote by
$S_m(x_1,...,x_k)$ the $m$-th elementary symmetric polynomial in $k$ commuting variables, such that
\bea
&& S_0(x_1,...,x_k)=1 \nonumber \\
&& S_1(x_1,...,x_k)=x_1+\cdots + x_k \nonumber \\
&& S_2(x_1,...,x_k)=\sum_{i<j} x_i x_j \nonumber \\
&& \cdots \nonumber \\
&& S_k(x_1,...,x_k)=x_1x_2 \cdots x_k. \nonumber
\eea
Throughout this section $p$ is odd.

We let
$$F_{2n+1,i,p}(t)={\rm Res}_{x_1,...,x_{2n+1}} \frac{1}{(x_1 \cdots x_{2n+1})^i} \prod_{i<j} (x_i-x_j)^p  \prod_{i=1}^{2n+1} (1+x_i)^t  \prod_{i=1}^n {\rm ln}(1-x_{2i}/x_{2i-1}),$$
with an eye on the special case $n=2$ and $i=5p$.

One can in principle try something more general, but we found that
$${\rm Res}_{x_1,...,x_{2n+1}} \frac{ x_{1}^{i_1} \cdots x_{2n+1}^{i_n} }{ (x_1 \cdots x_{2n+1})^i} \prod_{i<j} (x_i-x_j)^p  \prod_{i=1}^{2n+1} (1+x_i)^t  \prod_{i=1}^n {\rm  ln}(1-x_{2i}/x_{2i-1})$$
has rather nasty properties, unless of course $i_1=i_2=...=i_{2n+1}$ (which then reduces to $F_{2n+1,i-i_1,p}(t)$). This fact alone makes evaluation of $I_{2n+1,i,p}(t)$ difficult
by using the first order recursions where we go from $i_1$ to $i_i+1$, etc.
To fix this problem we introduce slightly more general residues. For $m \in \{0,...,2n+1\}$, and $S_m$ as above, we let
$$I_{2n+1,i,p,m}(t)$$
$$={\rm Res}_{x_1,...,x_{2n+1}} \frac{S_m(x_1,...,x_{2n+1})}{(x_1 \cdots x_{2n+1})^i} \prod_{i<j} (x_i-x_j)^p  \prod_{i=1}^{2n+1} (1+x_i)^t  \prod_{i=1}^n {\rm ln}(1-x_{2i}/x_{2i-1}).$$
These expressions are clearly polynomials in $t$. Actually, we can say more

\begin{lemma} \label{degree}
Let $i \in \mathbb{Z}$. For $i > pn+1$,
$${\rm deg}(I_{2n+1,i,p,m}(t)) \leq (i-pn-1)(2n+1)-m.$$
For $i=pn+1$,
$$I_{2n+1,pn+1,p,0}(t)=const$$
For $i <pn$,
$$I_{2n+1,i,p,m}(t)=0$$
\end{lemma}

The lemma (or easy inspection) gives
$${\rm deg}(I_{2n+1,i,p,m-1}(t))+1 \geq {\rm deg}(I_{2n+1,i,p,m}(t)),$$
but as we shall see we have strong equality here.
Also,
\be \label{link}
I_{2n+1,i,p,2n+1}(t)=I_{2n+1,i-1,p,0}(t),
\ee
which simply follows from
$$S_{2n+1}(x_1,...,x_{2n+1})=x_1 \cdots x_{2n+1}.$$
We want to focus on $m \in \{1,...,2n+1 \}$.
The goal is  to  obtain a recursion
\be \label{recursion}
I_{2n+1,i,p,m}(t)=(a_m t+b_m)I_{2n+1,i,p,m+1}(t),
\ee
where $a_m$ and $b_m$ depend on $n$,$p$ and $i$ as well.
Suppose we are able to do that. Then,   by iterating the relation (\ref{recursion}) we obtain
$$I_{2n+1,i,p,0}(t)=(a_1 t+b_1)I_{2n+1,i,p,1}(t)=\cdots=\prod_{i=1}^m(a_it+b_i)I_{2n+1,i,p,m}(t)=\cdots$$
$$=\prod_{i=1}^{2n+1}(a_it+b_i)I_{2n+1,i-1,p,2n+1}(t)=\prod_{i=1}^{2n+1}(a_it+b_i)I_{2n+1,i-1,p,0}(t).$$
Then we can continue iterating $I_{2n+1,i-1,p,0}(t)$, etc. This process will eventually terminate after finitely many steps  when we  reach $I_{2n+1,pn,p,0}(t)=const$
(see Lemma \ref{degree}).

In order to get  (\ref{recursion}) we need several auxiliary results. We will be using
$$S_m(x_1,...,\hat{x_i},...,x_{2n+1})$$
to denote the $m$-th elementary symmetric function in all variables except $x_i$. Similarly, we denote
by
$$S_m(x_1,...,\hat{x_i},..,\hat{x_j},..,x_{2n+1}) $$
$m$-th elementary symmetric function in all variables except $x_i$ and $x_j$.
For instance
$$S_3(x_1,x_2,\hat{x_3},x_4,x_5)=x_1 x_2 x_4+x_1 x_2 x_5+x_1 x_4 x_5+x_2x_4x_5.$$

Next couple of elementary lemmas can be easily proven, so we omit the proofs.
\begin{lemma} We have
\bea
&& \sum_{i=1}^{2n+1} x_i S_m(...,\hat{x_i},...)=(m+1)S_{m+1}(x_1,...,x_{2n+1}) \nonumber \\
&& \sum_{i=1}^{2n+1} S_m(...,\hat{x_i},...)=(2n+1-m)S_{m}(x_1,...,x_{2n+1}) \nonumber 
\eea
\end{lemma}

Similarly,
\begin{lemma} Assume $m \geq 2$. Then:
\bea
&& \sum_{i<j} x_i x_j S_{m-1}(...,\hat{x_i},..,\hat{x_j},...)=  \frac{m(m+1)}{2} S_{m+1}(x_1,...,x_{2n+1}) \nonumber \\
&& \sum_{i <j} S_m(...,\hat{x_i},..,\hat{x_j},...)= \frac{(2n-m+1)(2n-m)}{2}  S_{m}(x_1,...,x_{2n+1}) \nonumber \\
&& \sum_{i<j} (x_i+x_j)S(...,\hat{x_i},..,\hat{x_j},...)=(m+1)(2n-m)S_{m+1}(x_1,...,x_{2n+1}) \nonumber
\eea
\end{lemma}

\begin{lemma} \label{Jm} For $i=1,...,n$, we have
$$\left((1+x_{2i})x_{2i} S_m(...,\hat{x_{2i}},...) \frac{\partial}{\partial  x_{2i}}+(1+x_{2i-1})x_{2i-1} S_m(...,\hat{x_{2i-1}},...)\frac{\partial}{\partial  x_{2i-1}} \right) {\rm ln}(1-x_{2i}/x_{2i-1})$$
$$=x_{2i}(S_{m}(...,\hat{x_{2i-1}},\hat{x_{2i}},...)-S_{m-1}(...,\hat{x_{2i-1}},\hat{x_{2i}},...)):=J_{i,m} $$

\end{lemma}

\begin{example} For $m=1$ and $n=1$ and $i=1$:
$$J_{1,1}=x_2(x_3-1).$$
For $m=1$ and $n=2$ and $i=1$,
$$J_{1,1}=x_2(x_3x_4+x_3x_5+x_4x_5-x_3-x_4-x_5).$$
\end{example}

\begin{lemma}

$$\left\{(1+x_i)x_i S_m(...,\hat{x_i},...)\frac{\partial}{\partial x_i}+(1+x_j)x_j S_m(...,\hat{x_j},...)\frac{\partial}{\partial x_j} \right\} (x_i-x_j)^p $$
$$=p\left\{x_i x_j S_{m-1}(...,\hat{x_i},..,\hat{x_j},...)+S_m(...,\hat{x_i},..,\hat{x_j},...)+ (x_i+x_j)S(...,\hat{x_i},..,\hat{x_j},...)\right\}(x_i-x_j)^p.$$

\end{lemma}

\begin{lemma}
For every $i \in \{1,...,n \}$, we have
\be \label{res0}
{\rm Res}_{x_1,...,x_{2n+1}} \frac{1}{(x_1 \cdots x_{2n+1})^i} \prod_{i<j} (x_i-x_j)^p  \prod_{i=1}^{2n+1} (1+x_i)^t J_{i,m}  \prod_{j=1, j \neq i }^n {\rm ln}(1-x_{2j}/x_{2j-1})=0.
\ee

\end{lemma}

\begin{proof} Let us take $i=1$ for simplicity and $m$ at least two.  Similar argument holds in general.  From Lemma \ref{Jm} we see that
$$J_{1,m}=x_{2} \left\{S_{m}(\hat{x_{1}},\hat{x_{2}},...)-S_{m-1}(\hat{x_{1}},\hat{x_{2}},...)\right\}$$ does not involve $x_1$ variable.
We show that $J_{i,m}$ in (\ref{res0}) can be replaced with either $x_{2}S_{m}(\hat{x_{1}},\hat{x_{2}},...)$ of $x_{2}S_{m-1}(\hat{x_{1}},\hat{x_{2}},...)$ and the result is still zero.
To see that observe that each expression is a sum of $m$-th (or $m-1$-th) elementary monomials not involving variable $x_1$ and $x_2$. Such a monomial either
Case 1. involve the last variable $x_{2n+1}$, or Case 2. has no variable $x_{2n+1}$. In the first case we apply $x_1 \leftrightarrow x_{2n+1}$ and Lemma \ref{trivial} below.
In Case 2 we use $x_2 \leftrightarrow x_{2n+1}$ and the lemma below. Observe that logarithmic expression
$$ \prod_{ j \neq 1 }^n {\rm ln}(1-x_{2j}/x_{2j-1})$$
is invariant under this permutation (it does not involve $x_1$, $x_2$ or $x_{2n+1}$).
\end{proof}

\begin{lemma} \label{trivial}
$${\rm Res}_x  \frac{d}{dx} (f(x))=0.$$
$${\rm Res}_{x_1,...,x_{2n+1}} \frac{1}{(x_1 \cdots x_{2n+1})^i} \prod_{i<j} (x_i-x_j)^p=0$$
\end{lemma}
\begin{proof}The first formula is clear. For the second relation observe that the substitution $x_i \leftrightarrow x_j$ adds an additional sign
so we have
$${\rm Res}_{x_1,...,x_{2n+1}} \frac{1}{(x_1 \cdots x_{2n+1})^i} \prod_{i<j} (x_i-x_j)^p=-{\rm Res}_{x_1,...,x_{2n+1}} \frac{1}{(x_1 \cdots x_{2n+1})^i} \prod_{i<j} (x_i-x_j)^p.$$
\end{proof}

\begin{theorem} \label{main-thm}Let $m \in \{0,...,2n \}$ and $n \geq 1$. Then we have
\bea \label{main}
&& \frac{1}{2}(m-2n-1)(pm+2i-2-2pn) I_{2n+1,i,p,m}(t)  \\
&& =(m+1)(t+p(2n-m)+\frac{pm}{2}-i+2)I_{2n+1,i,p,m+1}(t) \nonumber
\eea

\end{theorem}
\begin{proof}  Observe that
$${\rm Res}_{x_1,...,x_{2n+1}} \sum_{i=1}^{2n+1} \frac{\partial}{\partial x_{i}}  \biggl( (1+x_i)x_i S_m(...,\hat{x_i},...) \frac{1}{(x_1 \cdots x_{2n+1})^i} \prod_{i<j} (x_i-x_j)^p $$
$$\prod_{i}(1+x_i)^t \prod_{i=1}^n {\rm ln}(1-x_{2i}/x_{2i-1}) \biggr)=0$$

Now we apply the product rule for differentiation and the above lemmas  to simplify the resulting expression down to

\bea 
&& (2n+1-m)I_{2n+1,i,p,m}(t)+2(m+1)I_{2n+1,i,p,m+1}(t)-i(2n+1-m)I_{2n+1,i,p,m}(t) \nonumber \\
&& -i(m+1)I_{2n+1,i,p,m+1}(t)+\frac{p(2n-m+1)(2n-m)}{2} I_{2n+1,i,p,m}(t) \nonumber \\
&& +\frac{pm(m+1)}{2} I_{2n+1,i,p,m+1}(t)+p(m+1)(2n-m)I_{2n+1,i,p,m+1}(t)+t(m+1)I_{2n+1,i,p,m+1}(t)=0. \nonumber
\eea
When we simplify this we obtain the wanted formula.
\end{proof}

\begin{corollary} \label{cor-5}
Conjecture 10.1 in \cite{AdM-IMRN} holds, i.e.,
$$ I_{5, 5p, p,0} (t) =  B_p {t \choose 3p-1} {t+p/2 \choose 3p-1} {t+p \choose 3p-1} {t+3p/2 \choose 3p-1} {t+2p \choose 3p-1} \quad (B_p \ne 0).$$
\end{corollary}

\begin{proof}  Take $n=2$ and $i=5p$. An application of Theorem \ref{main-thm} with $m=0,1,2,3,4$ gives a polynomial of degree $5$. This reduces
$i$ down to $5p-1$. Repeat this procedure until we reach $i=2p+1$. The resulting $t$ polynomial is of degree $5(3p-1)=15p-5$, and it is easy to see that this polynomial is (up to a constant) proportional
to
$${t \choose 3p-1} {t+p/2 \choose 3p-1} {t+p \choose 3p-1} {t+3p/2 \choose 3p-1} {t+2p \choose 3p-1}.$$
Non-vanishing of  $B_p$ was proved in  Theorem 6.3 and Proposition 6.4  of \cite{AdM-IMRN}.

\end{proof}

\begin{corollary} Theorem \ref{theorem-ct} holds.
\end{corollary}
\begin{proof} Follows along the lines of Corollary \ref{cor-5}. 
Non-vanishing of $\l_{p,k}$ is again consequence of  Theorem 6.3  of \cite{AdM-IMRN}
and the fact that  $G^k  Q e ^{ - (2 k+1) \alpha}$ is a non-trivial singular 
vector.
\end{proof}

\begin{corollary}
Inside the Zhu algebra $A(\tripletwp)$ we have
\bea f_{2,p} (x) &:=&\left( \prod _{i =1} ^{3p-1} (x - h  _{1,i} )  \right) \left(  \prod_{i =1} ^{
 \tfrac{3p-1}{2}
} (x - h _{1,2p-i} )\right) \left(\prod_{i =1} ^{3p-1} (x - h  _{2,i} )\right )=0, \nonumber
\eea
where $x=[\omega]$.
\end{corollary}

Regarding the constant term (which would allow us to explicitly determined the nonzero constant) we
can only conjecture (as in \cite{AdM-IMRN}):

\begin{conjecture} \label{log-dyson} Let $k$ and $m$ be positive integers. Then (up to a sign)
$${\rm CT}_{x_1,...,x_{2k+1}}  \frac{1}{(x_1 \cdots x_{2k+1})^{(2m+1)k}} \prod_{i=1}^k {\rm ln}\left(1-\frac{x_{2i}}{x_{2i-1}}\right)  {\prod_{1 \leq i <j \leq 2k+1
}(x_i-x_j)^{2m+1}}$$ equals
$$\frac{((2k+1)(2m+1))!!}{(2k+1)!!(2m+1)!!^{2k+1}}.$$
\end{conjecture}

\noindent This conjectural formula is what we call {\em logarithmic Dyson's constant term identity}.
\begin{remark} We were able to settle the conjecture for $m=1$ for all $p$ \cite{AdM-IMRN}, and for $m=2$ (also for all $p$) by using different methods.
\end{remark}

\section{Graded dimensions}

In this part we analyze irreducible $\WW_{2,p}$ graded dimensions (or simply, characters) 
and their $SL(2,\mathbb{Z})$-closure.
We introduce some notation first.

For an $L(0)$-diagonalizable module $M$ with finite-dimensional eigenspaces we define its character
$$\chi_M(q)={\rm tr}_{M} q^{L(0)-c/24}, \ \ q=e^{2 \pi i \tau},  \ \ | q| <1.$$
We also recall the usual Jacobi theta constants 
$$\theta_{r,2p} (\tau)=\sum_{j \in \frac{r}{4p}+\mathbb{Z}} q^{2pj^2}.$$
For simplicity we let $\theta_{r}(\tau):=\theta_{r,2p} (\tau)$. 
Derivatives of theta constants  will be denoted by $\theta'$, $\theta''$, etc.

It is easy to see that $$\theta_{r+4p}(\tau)=\theta_{r}(\tau)=\theta_{-r,p}(\tau),$$ and
 $$\theta'_{r}(\tau)=-\theta'_{-r}(\tau), \ \ \theta''_{r}(\tau)=\theta''_{-r}(\tau).$$
This implies, $$\theta'_0(\tau)=\theta'_{2p}(\tau)=0.$$
The following result (proven in \cite{FGST-log}) follows from 
decomposition of irreducible $\WW_{2,p}$-modules as Virasoro representations.

\begin{proposition}
Denote by $\chi_{r,s}^\pm(q)$ the character of $\X_{r,s}^\pm$, and by $\chi_{r,s}(q)$ 
the character of $L(c_{2,p},h_{r,s})$. Then we have 
\bea \label{char-vir}
&& \chi_{r,1}(q)=\frac{1}{\eta}(\theta_{p-2r}(\tau)-\theta_{p+2r}(\tau)), \ \ 1 \leq r \leq \frac{p-1}{2} \nonumber \\
&& \chi^+_{r,s}(q)=\frac{1}{4p^2 \eta}(\theta''_{p-2r}(\tau)-\theta''_{p+2r}(\tau)
-(ps+2r)\theta'_{ps+2r}(\tau)+(ps-2r)\theta'_{ps+2r}(\tau)  \nonumber \\
&& + \frac{(ps+2r)^2}{4}\theta_{ps+2r}(\tau)-\frac{(ps-2r)^2}{4}\theta_{ps-2r}(\tau)), \nonumber \\
&& \chi^-_{r,s}(q)=\frac{1}{4p^2 \eta}(\theta''_{2p-ps-2r}(\tau)-\theta''_{2p+ps-2r}(\tau)
+(ps+2r)\theta'_{2p-ps-2r}(\tau)+(ps-2r)\theta'_{2p+ps-2r}(\tau) \nonumber \\
&&  + \frac{(ps+2r)^2-4p^2}{4}\theta_{2p-ps-2r}(\tau)-\frac{(ps-2r)^2-4p^2}{4}\theta_{2p+ps-2r}(\tau)), 
\nonumber \eea
where $1 \leq r \leq p, 1 \leq s \leq 2$.
\end{proposition}

The vector space spanned by $\chi^\pm_{r,s}$ and $\chi_{r,1}$ is obviously too small to be modular invariant. Instead we have
\begin{theorem} \label{mod-basis}
The $SL(2,\mathbb{Z})$-closure of the space of irreducible $\WW_{2,p}$ characters is $\frac{15p-5}{2}$-dimensional with a basis:
\bea
&& \mathcal{B}= \mathcal{B}_2 \cup \mathcal{B}_1 \cup \mathcal{B}_0 \nonumber \\
&& \mathcal{B}_2=\left\{ \frac{\tau^i (\theta''_{p-r}(\tau)-\theta''_{p+r}(\tau))}{\eta}, 1 \leq r \leq \frac{p-1}{2}, 0 \leq i \leq 2 \right\} \nonumber \\
&& \mathcal{B}_1=\left\{ \frac{\tau^i \theta'_{r}(\tau)}{\eta},  \ \ 1 \leq r \leq 2p-1, \  0 \leq i \leq 1 \right\}, \ \ \mathcal{B}_0= \left\{ \frac{\theta_{r}(\tau)}{\eta}, \ \ 0 \leq r \leq 2p  \right\}. \nonumber
\eea
\end{theorem}

\begin{proof} Let us denote the vector space spanned by $\mathcal{B}=\mathcal{B}_2 \cup \mathcal{B}_1 \cup \mathcal{B}_0$ with
 $\mathcal{F}_p$. It is easy to see ${\rm dim}(\mathcal{F}_p)=\frac{3(p-1)}{2}+2(2p-1)+2p+1=\frac{15p-5}{2}$. 
To show modular invariance of $\mathcal{F}_p$ it is sufficient to observe that each $\mathcal{B}_i \subset \mathcal{B}$
generates a modular invariant subspace, which is easy to check directly.
To finish the proof we must show that each element in $\mathcal{B}$ can be obtained by applying appropriate 
modular transformation (on a linear combination) of irreducible characters. 
Since characters of irreducible $V_L$-modules are given by theta quotients $\frac{\theta_r}{\eta}$, $0 \leq r \leq 2p$, and those
are linear combinations of irreducible $\WW_{2,p}$ characters, every element in 
$\mathcal{B}_0$ is contained in the closure. Similarly, by taking linear combinations of $\X^+_{r,s}$ and $\X^-_{r,s}$ we 
see that $\frac{\theta'_{r}}{\eta}$ is also contained in the closure. But then the relation
$$\frac{\theta'_{r}}{\eta}|_{\tau \mapsto -\frac{1}{\tau}}=\tau \sum_{r=0}^{2p} \alpha_k \frac{\theta'_k}{\eta},$$
where $\alpha_k$ are some constants, implies $\mathcal{B}_1$ is also in the closure. 
Now, $f_r:=\frac{\theta''_{p-r}(\tau)-\theta''_{p+r}(\tau)}{\eta}$ are also in the closure and again by applying $\tau \mapsto \frac{-1}{\tau}$ 
on $f_r$ we see that all of $\mathcal{B}_2$ are in fact in the closure. The proof follows.
\end{proof}

\begin{remark} In \cite{FGST-log} a similar (but different) basis of the $SL(2,\mathbb{Z})$ closure was constructed, motivated 
by description of the center of $\mathcal{U}_{2,p}$.
\end{remark}


\begin{thebibliography}{FGST3}


\bibitem[A]{A-2003} D. Adamovi\'{c},
Classification of irreducible modules of certain subalgebras of
free boson vertex algebra, J. Algebra 270 (2003) 115-132.


\bibitem[AM1]{AdM-2007} D. Adamovi\'c and A. Milas, Logarithmic intertwining operators
and $\mathcal{W}(2,2p-1)$-algebras, {\em Journal of Math. Physics}
{\bf 48}, 073503 (2007).

\bibitem[AM2]{AdM-triplet} D. Adamovi\'c and A. Milas, On the triplet vertex algebra
$\mathcal{W}(p)$, {\em Advances in Math.} 217 (2008) 2664-2699;
arxiv:0707.1857.

\bibitem[AM3]{AdM-lattice} D. Adamovi\'c and A. Milas, {\em  Lattice
construction of logarithmic modules for certain vertex algebras}, 
{\em Selecta Mathematics (NS)},   {\bf 15} (2009) 535-561; {\tt arxiv}:0902.3417.

\bibitem[AM4]{AdM-IMRN}  D. Adamovi\'c and A. Milas, On $W$-algebras associated to $(2,p)$ minimal models and their representations,  {\em IMRN}, (2010) {\bf 20}, 3896-3934.

\bibitem[AM5]{AdM-2010}  D. Adamovi\'c and A. Milas, The structure of Zhu's algebras for certain $\W$-algebras, arXiv:1006.5134v1, submitted for publication.

\bibitem[AM6]{AdM-2011} D. Adamovi\'c and A. Milas, An explicit realization of logarithmic modules for triplet vertex algebra ${\mathcal W}_{p,p'}$, to appear

\bibitem[DL]{DL} C. Dong and J.Lepowsky,  
{\em Generalized vertex algebras and relative vertex operators}, Birkh\"auser,  Boston, 1993.



\bibitem [FGST1]{FGST-log} B.L. Feigin, A.M. Ga\u\i nutdinov, A. M. Semikhatov, and I. Yu Tipunin,
 Logarithmic extensions of minimal
models: characters and modular transformations, Nucl. Phys. B 757
(2006) 303�343.

\bibitem[FGST2]{FGST-log-kl}  B.L. Feigin, A.M. Ga\u\i nutdinov, A. M. Semikhatov, and I. Yu Tipunin, Kazhdan--Lusztig-dual quantum group for logarithmic extensions of Virasoro minimal models
,  {\em Jour. of Math. Physics } {\bf 48}  032303 (2007)

\bibitem[FW]{FW} P. Forrester,  O. Warnaar, The importance of the Selberg integral. {\em Bull. Amer. Math. Soc.}  (N.S.) {\bf 45} (2008), no. 4, 489--534.

\bibitem[GRW1]{GRW-1} M. Gaberdiel, I. Runkel and S. Wood, Fusion rules and boundary conditions in the $c=0$ triplet model, J. Phys. A: Math. Theor. 42 (2009) 325--403,
{\tt arxiv}: 0905.0916.

\bibitem[GRW2]{GRW-2} M. Gaberdiel, I. Runkel and S. Wood, A modular invariant bulk theory for the c=0 triplet model
, J.Phys. A: Math. Theor. 44 (2011) 015204
, {\tt arxiv}: 1008.0082v1

\bibitem[Hu]{Hu} Y.-Z. Huang, {\em Cofiniteness conditions, projective covers and the logarithmic tensor product theory},
J. Pure Appl. Algebra {\bf 213} (2009) 458-475;
 {\tt arxiv.0712.4109}.

\bibitem[HLZ]{HLZ} Y.-Z. Huang,  J. Lepowsky, and L. Zhang, {\em Logarithmic tensor product theory for generalized modules for a conformal vertex algebra}, {\tt
arXiv:0710.2687}.


\bibitem[LL]{LL} J. Lepowsky and H. Li, {\em Introduction to Vertex Operator
Algebras and Their Representations}, Progress in Mathematics, Vol.
227, Birkh\"auser, Boston, 2003.

\bibitem[M]{M} W.G. Morris, Constant term identities for finite and affine root systems, PhD thesis, UWisconsin, Madison, 1982.

\bibitem[Miy]{Miy} M. Miyamoto, Flatness of Tensor Products and Semi-Rigidity for $C_2$-cofinite Vertex Operator Algebras I and II, preprints; {\tt arXiv}:0906.1407 and {\tt arXiv}: 0909.3665.

\bibitem[PRZ]{PRZ}  P. Pierce, J. Rasmussen and J.-B. Zuber, Logarithmic minimal models, {\em J.Stat.Mech.} {\bf 611} :P017,2006.

\bibitem[R]{R} J. Rasumussen, W-extended logarithmic minimal models, {\em Nucl.Phys.} B {\bf 807}, (2009), 495-533.

\bibitem[TK]{TK} A. Tsuchiya and Y. Kanie, Fock space representations of
the Virasoro algebra - Intertwining operators,  { Publ. RIMS} {\bf
22} (1986), 259-327.

\bibitem[W]{W} S. Wood,  Fusion Rules of the $\mathcal{W}_{p,q}$ Triplet Models, J. Phys. A: Math. Theor. 43 (2010) 045212, {\tt  arXiv:0907.4421}.

\end{thebibliography}
\end{document}